\documentclass[]{amsart}
\usepackage{amscd,amsthm,amssymb,amsfonts,amsmath,euscript}
\theoremstyle{plain}
\newtheorem{thm}{Theorem}[section]

\newtheorem{prop}[thm]{Proposition}
\newtheorem{cor}[thm]{Corollary}

\theoremstyle{definition}

\theoremstyle{remark}
\newtheorem{rem}[thm]{Remark}
\newtheorem*{thank}{Acknowledgments}
\newcommand{\nc}{\newcommand}

\def\makeop#1{\expandafter\def\csname#1\endcsname
  {\mathop{\rm #1}\nolimits}\ignorespaces}
\makeop{Hom}   \makeop{End}   \makeop{Aut}   \makeop{Isom}  \makeop{Pic} 
\makeop{Gal}   \makeop{ord}   \makeop{Char}  \makeop{Div}   \makeop{Lie} 
\makeop{PGL}   \makeop{Corr}  \makeop{PSL}   \makeop{sgn}   \makeop{Spf}
\makeop{Spec}  \makeop{Tr}    \makeop{Nr}    \makeop{Fr}    \makeop{disc}
\makeop{Proj}  \makeop{supp}  \makeop{ker}   \makeop{im}    \makeop{dom}
\makeop{coker} \makeop{Stab}  \makeop{SO}    \makeop{SL}    \makeop{SL}
\makeop{Cl}    \makeop{cond}  \makeop{Br}    \makeop{inv}   \makeop{rank}
\makeop{id}    \makeop{Fil}   \makeop{Frac}  \makeop{GL}    \makeop{SU}
\makeop{Trd}   \makeop{Sp}    \makeop{Tr}    \makeop{Trd}   \makeop{diag}
\makeop{Res}   \makeop{ind}   \makeop{depth} \makeop{Tr}    \makeop{st}
\makeop{Ad}    \makeop{Int}   \makeop{tr}    \makeop{Sym}   \makeop{can}
\makeop{length}\makeop{SO}    \makeop{torsion} \makeop{GSp}
\def\makebb#1{\expandafter\def
  \csname bb#1\endcsname{{\mathbb{#1}}}\ignorespaces}
\def\makebf#1{\expandafter\def\csname bf#1\endcsname{{\bf
      #1}}\ignorespaces} 
\def\makegr#1{\expandafter\def
  \csname gr#1\endcsname{{\mathfrak{#1}}}\ignorespaces}
\def\makescr#1{\expandafter\def
  \csname scr#1\endcsname{{\EuScript{#1}}}\ignorespaces}
\def\makecal#1{\expandafter\def\csname cal#1\endcsname{{\mathcal
      #1}}\ignorespaces} 
\def\doLetters#1{#1A #1B #1C #1D #1E #1F #1G #1H #1I #1J #1K #1L #1M
                 #1N #1O #1P #1Q #1R #1S #1T #1U #1V #1W #1X #1Y #1Z}
\def\doletters#1{#1a #1b #1c #1d #1e #1f #1g #1h #1i #1j #1k #1l #1m
                 #1n #1o #1p #1q #1r #1s #1t #1u #1v #1w #1x #1y #1z}
\doLetters\makebb   \doLetters\makecal  \doLetters\makebf
\doLetters\makescr 
\doletters\makebf   \doLetters\makegr   \doletters\makegr
     \def\qed{\qedmark\medbreak}%
\def\qedmark{{\enspace\vrule height 6pt width 5pt depth 1.5pt}}%

\normalsize

\makeop{Bl}

\def\Qbar{\overline{\bbQ}}

\newcommand{\Z}{\mathbb Z}

\newcommand{\C}{\mathbb C}
\newcommand{\F}{\mathbb F}

\nc{\embed}{\hookrightarrow}
\newcommand{\ch}{characteristic }

\nc{\ol}{\overline}
\nc{\wt}{\widetilde}
\nc{\opp}{\mathrm{opp}}
\nc{\ul}{\underline}

\makeop{Ram}
\makeop{Rep}

\begin{document}
\title{Irreducibility of the Siegel moduli spaces with parahoric level
  structure} 
\author{Chia-Fu Yu}
\address{
Institute of Mathematics \\
Academia Sincia \\
128 Academia Rd.~Sec.~2, Nankang, Taipei, Taiwan }
\email{chiafu@math.sinica.edu.tw}
\date{\today}

\begin{abstract}
  We prove that the moduli space $\calA_{g,\Gamma_0(p)}\otimes \ol
  \F_p$ of principally polarized abelian varieties of dimension $g$ with a
  $\Gamma_0(p)$-level structure in \ch $p$ has $2^g$ irreducible components.
  The cases of other parahoric level structures are also considered.
\end{abstract} 


\maketitle

%\tableofcontents   % Table of Contents
% {\small
% 2000 MSC: Primary: 11G18; Secondary: 14G35. 

% Keywords and phrases: Siegel moduli schemes, parahoric level
% structure, $p$-adic monodromy. }

\section{Introduction}
\label{sec:01}
Let $p$ be a rational prime number. Let $N\ge 3$ be a prime-to-$p$
positive integer. We choose a primitive $N$th root of unity $\zeta_N$ in
$\Qbar\subset \C$ and an embedding $\Qbar\subset \Qbar_p$. Let
$\calA_{g,1,N}$ denote the moduli space over $\Z_{(p)}[\zeta_N]$ of
principally polarized abelian varieties with a full symplectic level $N$
structure with respect to $\zeta_N$. The moduli scheme $\calA_{g,1,N}$
has irreducible geometric fibers. A principally polarized abelian
scheme over a base scheme $S$ with a $\Gamma_0(p)$-level structure
consists of  
\[ (A,\lambda, 0=H_0\subset H_1\subset H_2\subset\dots \subset H_g\subset
A[p]\,),\]
where $(A,\lambda)$ is a principally polarized abelian scheme, $H_i$
is a locally free finite subgroup scheme of rank $p^i$ and $H_g$ is
isotropic for the Weil pairing $e_{\lambda}$ induced by $\lambda$. 
Let $\calA_{g,\Gamma_0(p),N}$ denote the moduli space over
$\Z_{(p)}[\zeta_N]$ that 
  parameterizes the objects in $\calA_{g,1,N}$ with a
  $\Gamma_0(p)$-level structure. We prove 
\begin{thm}\label{11}
The reduction $\calA_{g,\Gamma_0(p),N}\otimes \ol \F_p$ of the moduli
space modulo $p$ has $2^g$ irreducible components.    
\end{thm}
The case of $g=1$ is a well-known result of Deligne and Rapoport
\cite{deligne-rapoport}, where two irreducible components intersect at
supersingular points transversally. The case of $g=2$ has been proved by
de Jong \cite{dejong:gamma}. He also 
proves the flatness of the structure morphism and determines the
singularities. G\"ortz further proves the flatness and determines the
singularities for arbitrary genus $g$ (see \cite[Theorem
2.1]{goertz:symplectic} for precise statements). Particularly he
verifies the conjecture of Rapoport and Zink 
\cite{rapoport-zink} for the case of unramified symplectic groups.    

In \cite{ngo-genestier:alcoves} Genestier and Ng\^o study the
stratification by $\mu$-permissible elements in the sense of
Kottwitz-Rapoport on the moduli space (see
\cite{kottwitz-rapoport:alcoves}).  
They show that the $p$-rank is
constant on each stratum and determine the dimensions of the
strata. Particularly they prove
\begin{thm} [Genestier-Ng\^o] \label{12}
The ordinary locus $\calA^{\rm{ord}}$ of
$\calA:=\calA_{g,\Gamma_0(p),N}\otimes \ol \F_p$ is dense in
$\calA$. Moreover $\calA^{\rm{ord}}$ is the smooth locus of $\calA$.
\end{thm}

More precisely, let 
\[ \calA=\coprod_{w} \calA_w \]
be the stratification of $\calA$ by $\mu$-permissible elements.
Genestier and Ng\^o show that each stratum $\calA_w$ is non-empty and
that $p$-${\rm rank}(w)=g$ if and only if 
$\dim \calA_w=\dim \calA$. There are $2^g$ of such $\mu$-permissible
elements. The task taken here then is to show that each maximal
stratum $\calA_w$ is irreducible.
Indeed we give the simple combinatorial
description for such $\mu$-permissible elements by subsets $\tau$ of
$\{1,\dots, g\}$ and then show that the corresponding stratum
$\calA_\tau$ (see the definition in the next section) is irreducible. 
Our basic tool is the $p$-adic monodromy developed in
\cite{faltings-chai}.

\section{Proof of the theorem}
\label{sec:02}

We keep the notation as in the previous section. Let $S$ be a locally
noetherian base scheme in \ch $p$. Set
$H:=(\Z/p\Z\times \mu_p)^g$ and let $\varphi:H\times H\to \mu_p$ be the
alternating pairing defined by 
\[ \varphi( (m_i,\zeta_i),(n_i,\eta_i))=\prod_{i=1}^g \eta_i^{m_i}
\zeta_i^{-n_i}, \quad \forall\,m_i,n_i\in \Z/p\Z, \ \zeta_i,\eta_i\in
\mu_p. \]  
A $\Gamma(p)$-level structure on an ordinary principally polarized
abelian scheme $(A,\lambda)$ over $S$ is an isomorphism $\xi:H_S\simeq
A[p]$ such that the pull-back by $\xi$ of the Weil pairing
$e_\lambda$ is $\varphi$. 

Let $\calA^{\rm ord}_{\Gamma(p)}$ denote the cover of 
$\calA^{\rm ord}_{g,1,N}$ where the objects are endowed with a
$\Gamma(p)$-level structure. As a $\Gamma(p)$-level structure $\xi$
is determined by its restriction on $(\Z/p\Z)^g$,
this moduli scheme $\calA^{\rm ord}_{\Gamma(p)}$ represents the sheaf 
\[ \ul{Isom}((\Z/p\Z)^g,\calX[p]^{\text{\'et}}) \]
on $\calA^{\rm ord}_{g,1,N}$ for the \'etale topology, where
$(\calX,\lambda,\alpha)\to \calA^{\rm ord}_{g,1,N} $ is the universal
family. This cover gives rise to the map $\ol \rho:\pi_1(\calA^{\rm
  ord}_{g,1,N},\bar x)\to \Aut (\calX_{\bar x}[p]^{\text{\'et}})
  =\GL_g(\Z/p\Z)$.
It follows form the surjectivity of the 
$p$-adic monodromy $\rho:\pi_1(\calA^{\rm ord}_{g,1,N},\bar x)\to \GL_g(\Z_p)$
\cite{faltings-chai} that the moduli scheme 
$\calA^{\rm ord}_{\Gamma(p)}$ is irreducible.

Let $(A,\lambda, H_\bullet,\alpha)$ be an object in $\calA^{\ord}$ over a
{\it connected} scheme $S$. We associate a subset $\tau$ of
$\{1,\dots, g\}$ as follows: for $1\le i\le g$, $i\in \tau$ if and
only if $H_i/H_{i-1}$ is of multiplicative type. Therefore, the
ordinary locus 
\[ \calA^{\rm ord}=\coprod_{\tau} \calA_\tau \]
is a finite union of open and closed subschemes $\calA_\tau$ whose objects have
type $\tau$. 

\begin{prop}\label{21}
  For each $\tau\subset \{1,\dots,g\}$, there is a natural finite surjective
  morphism from $\calA^{\rm ord}_{\Gamma(p)}$ to $\calA_\tau$. 
\end{prop}
\begin{proof}
Given any $\tau$, we define a subgroup scheme $H_\tau:=\prod_{i=1}^g
K_i$ of $H$, where 
\[ K_i:=
\begin{cases}
  \Z/p\Z & \text{if $i\not\in \tau$,} \\
  \mu_p & \text{if $i\in \tau$.}
\end{cases} \]
It comes with a natural filtration of subgroup schemes
$H_{\tau,i}:=\prod_{j=1}^{i} K_j$. 

  Given an object $(A,\lambda,\xi,\alpha)_S$ of $\calA^{\rm
  ord}_{\Gamma(p)}$, let $\xi_\tau$ be the restriction of $\xi$ on
  $H_\tau\otimes S$. Note that $\xi$ is determined by its restriction
  $\xi_\tau$ by the Cartier duality. Put
  $H_i:=\xi_\tau(H_{\tau,i}\otimes S)$
  for $1\le i\le g$
  and then $H_\bullet$ gives a $\Gamma_0(p)$-level structure of type
  $\tau$ on $(A,\lambda)$. The map $(A,\lambda,\xi,\alpha)\mapsto
  (A,\lambda,H_\bullet,\alpha)$ defines a morphism $f_\tau$ from $\calA^{\rm
  ord}_{\Gamma(p)}$ to $\calA_\tau$. It is clear from the construction
  that $f_\tau$ is finite and surjective. \qed
\end{proof}

\begin{cor}\label{22}
  Every stratum $\calA_\tau$ is irreducible.
\end{cor}

Theorem~\ref{11} now is proved.

\section{Parahoric level structures}
\label{sec:03}
We keep the notation as before.
Let $\ul k=(k_1,\dots,k_r)$ be an ordered partition of $g$: each
$k_i \ge 1$ and $\sum_{j=1}^r k_j=g$. Set $h(i):=\sum_{j=1}^i k_j$ for
$1\le i\le r$. 
A $\Gamma_{\ul k}(p)$-level structure on a principally polarized abelian
scheme $(A,\lambda)$ is a flag of finite subgroup schemes
\[ 0=H_{h(0)}\subset H_{h(1)}\subset \dots\subset H_{h(r)}\subset
A[p], \]
where $H_{h(i)}$ is locally free of rank $p^{h(i)}$ and $H_{h(r)}$ is
isotropic for the Weil pairing $e_{\lambda}$. When $r=g$, it is a
$\Gamma_0(p)$-level structure defined in Section 1. 
Let $\calA_{g,\Gamma_{\ul k}(p),N}$ denote the moduli space over
$\Z_{(p)}[\zeta_N]$ that parameterizes the objects in $\calA_{g,1,N}$ with a
$\Gamma_{\ul k}(p)$-level structure. Let $\calA_{\ul k}:=
\calA_{g,\Gamma_{\ul k}(p),N} \otimes \ol \F_p$ be the reduction
modulo $p$, and let $\calA_{\ul k}^{\rm ord}$ denote the ordinary locus of
$\calA_{\ul k}$.
  
\begin{prop}\label{31}
  The ordinary locus $\calA_{\ul k}^{\rm ord}$ is open and dense in
  $\calA_{\ul k}$.
\end{prop}
\begin{proof}
  As the natural forgetful morphism $f:\calA\to \calA_{\ul k}$ is
  surjective, the statement follows from Theorem~\ref{12}. \qed 
\end{proof}

Let $(A,\lambda, H_\bullet,\alpha)$ be an object in $\calA^{\ord}_{\ul
  k}$ over a
{\it connected} scheme $S$. We associate an element $\ul
\tau=(\tau(1),\dots,\tau(r))$ of
$\prod_{i=1}^r I(k_i)$, where $I(k_i)=[0,k_i]\cap \Z$, as follows: 
for $1\le i\le r$ the local part
of  $H_{h(i)}/H_{h(i-1)}$ has rank $p^{\tau(i)}$. Therefore, the
ordinary locus 
\[ \calA^{\rm ord}_{\ul k}=\coprod_{\ul \tau} \calA_{\ul k,\ul \tau} \]
is a finite union of open and closed subschemes $\calA_{\ul k, \ul \tau}$
whose objects have type $\ul \tau$. 

\begin{prop}\label{32}
  For each $\ul \tau$ in $\prod_{i=1}^r I(k_i)$, there is a
  natural finite surjective morphism from $\calA^{\rm
  ord}_{\Gamma(p)}$ to $\calA_{\ul k,\ul \tau}$. Therefore, each
  stratum $\calA_{\ul k,\ul \tau}$ is non-empty and irreducible.
\end{prop}
\begin{proof}
  Given any element $\ul \tau$, we define a subgroup scheme $H_{\ul
  \tau}=\prod_{i=1}^r K_i$ of $H$, where $K_i=(\Z/p\Z)^{k_i-\tau(i)}\times
  (\mu_p)^{\tau(i)}$. It comes with a filtration of subgroup schemes
  $H_{\ul \tau,h(i)}=\prod_{j=1}^i K_i$. 

  Now given an object $(A,\lambda,\xi,\alpha)_S$ of $\calA^{\rm
  ord}_{\Gamma(p)}$, let $\xi_{\ul \tau}$ be the restriction of $\xi$ on
  $H_{\ul \tau}\otimes S$. Put $H_{h(i)}:=\xi_{\ul \tau}(H_{\ul
  \tau,h(i)}\otimes S)$ for $1\le i\le r$
  and then $H_\bullet$ gives a $\Gamma_{\ul k}(p)$-level structure of type
  $\ul \tau$ on $(A,\lambda)$. The map $(A,\lambda,\xi,\alpha)\mapsto
  (A,\lambda,H_\bullet,\alpha)$ defines a morphism $f_{\ul \tau}$ from
  $\calA^{\rm ord}_{\Gamma(p)}$ to $\calA_{\ul k,\ul \tau}$. It is
  clear again that the morphism $f_{\ul \tau}$ is finite and 
  surjective. \qed
\end{proof}

The ordinary locus $\calA_{\ul k}^{\rm ord}$ has $|\prod_{i=1}^r
I(k_i)|$ irreducible components. By Proposition~\ref{31}, we prove

\begin{thm}\label{33}
  The reduction $\calA_{\ul k}=\calA_{g,\Gamma_{\ul k}(p),N}\otimes
  \ol \F_p$ modulo $p$ has $\prod_{i=1}^r (k_i+1)$ irreducible components. 
\end{thm}

\begin{rem}
  Proposition~\ref{31} and Theorem~\ref{33} are obtained by Chai
  and Norman \cite[p.~273]{chai-norman:gamma2} when $r=1$.
\end{rem}
\begin{thank}
  The author thanks Genestier and Ng\^o for their inspiring work. He
  also wishes to thank the referee for careful reading 
  and helpful comments.
\end{thank}

\ \\

\begin{center}
  {\bf Erratum}
\end{center}
Page 2 line 21 is incorrect. The moduli scheme 
$\calA^{\rm  ord}_{\Gamma(p)}$ is not \'etale 
over $\calA^{\rm ord}_{g,1,N}$. However, $\calA^{\rm ord}_{\Gamma(p)}$
is still irreducible. The argument changes as follow. Let $\calA'$ be
a finite radical base change of $\calA^{\rm ord}_{g,1,N}$ such that
the extension
\[ 0\to \calX[p]^0\to \calX[p]\to \calX[p]^{\text{\'et}}\to 0 \]
splits. For example, take $\calA'$ to be $(\calA^{\rm
  ord}_{g,1,N})^{(p^{-n})}$ over $\F_q$. 
Let $\calX'$ be the base change of
  $\calX$ over $\calA'$ and $\calA'_{\Gamma(p)}$ be the finite cover
  of $\calA'$ classifying the $\Gamma(p)$-level structures. Then
  $\calA'_{\Gamma(p)}$ is \'etale over $\calA'$ and represents the sheaf
\[ \ul{Isom}((\Z/p\Z)^g,\calX'[p]^{\text{\'et}}) \]
on $\calA'$ for the \'etale topology. Then the same
argument shows $\calA'_{\Gamma(p)}$ irreducible. Since the natural
morphism $\calA'_{\Gamma(p)}\to \calA^{\rm  ord}_{\Gamma(p)}$ is
surjective, the target is irreduicble.

J.~Tilouine informed the author that $\calA^{\rm  ord}_{\Gamma(p)}$ is a
torsor over $\ul{Isom}((\Z/p\Z)^g,\calX[p]^{\text{\'et}})$ under the
finite group scheme $\mu_p\otimes_{\Z/p} U(\Z/p)$, where $U$ over
$\Z_p$ is the unipotent radical of the Siegel parabolic subgroup. We
thank him for helpful discussions. 

Page 3 line 7. Relace ``an ordered ... $g$'' by ``a tuple of positive
integers $k_i$ with $\sum_{j=1}^r k_j\le g$''. 
 
\end{document}